\documentclass[12pt]{article}

\usepackage{amssymb,bbm}
\usepackage{amsfonts}
\usepackage{latexsym}
\usepackage{amsmath}
\usepackage{amsthm}
\usepackage{geometry} 
\def\lddots{\mathinner{\mkern1mu\raise1pt\hbox{.}\mkern2mu  
\raise4pt\hbox{.}\mkern2mu\raise7pt\vbox{\kern7pt\hbox{.}}\mkern1mu}}
\makeatletter
\def\numberbysection{\@addtoreset{equation}{section}
 \def\theequation{\thesection.\arabic{equation}}}
 
\makeatother  
  
\numberbysection

\newcommand{\be}{\begin{eqnarray}}  
\newcommand{\ee}{\end{eqnarray}} 
\newcommand{\nn}{\nonumber}

\def\ds{\displaystyle}

\def\mf{\mathfrak}
\def\bb{\mathbbm}
\def\tbf{\textbf}
\def\C{\bb C}

\def\1{\bb{1}}

\def\pd{\partial}

\def\l{\lambda}

\def\g{\gamma}

\def\bp{\begin{proof}[{\bf Proof}]}

\newcounter{prop}
\newcounter{appendice}

\begin{document}

\begin{titlepage}



\vspace{5 mm}

\begin{center}
{\Large\bf Parametrization of semi-dynamical quantum reflection algebra}

\vspace{10 mm}

{\bf Jean Avan, Genevi\`eve Rollet}

\vspace{15mm}

{\it Laboratoire de Physique Th\'eorique et Mod\'elisation\\
Universit\'e de Cergy-Pontoise (CNRS UMR 8089), Saint-Martin 2\\
2  avenue Adolphe Chauvin, F-95302 Cergy-Pontoise Cedex, France.\\}
e-mail: {avan, rollet@u-cergy.fr}\\

\vspace{5mm}

\end{center}

\vspace{18mm}

\begin{abstract}
\noindent

We construct sets of structure
matrices for the semi-dynamical reflection algebra,
solving the Yang-Baxter type
consistency equations extended by
the action of an automorphism of the auxiliary space. These
solutions are 
parametrized by dynamical conjugation matrices,
Drinfel'd twist representations and quantum non-dynamical $R$-matrices.
They yield factorized forms for the monodromy matrices.
\end{abstract}

\vfill

\end{titlepage}
\setcounter{footnote}{0}
\section{Introduction}
The semi-dynamical reflection algebra (SDRA) was first formulated on a specific example in~\cite{ACF}.
The general formulation, together with a set of sufficient consistency conditions of Yang-Baxter type, was achieved in~\cite{NAR}.
The transfer matrix, commuting trace formulae, and representations of the
comodule structures, were defined in the same and in the following paper~\cite{NADR}; applications to the explicit construction of spin-chain type integrable Hamiltonians were given in~\cite{NA}.

The generators of the SDRA are encapsulated in a matrix $K(\l)$ acting on a vector space ${\cal V}$
denoted ``auxiliary space''. Two different types of auxiliary spaces will be considered
here: either a finite dimensional complex vector space $V$, or a loop space $V \otimes \C[[u]]$,
with $u$ the spectral parameter, in this last case the matrix $K(\l)$ should actually be denoted $K(\l,u)$ and belongs
to $\mbox{End}(V) \otimes \C[[u]]$. This matrix $K(\l)$ satisfies the semi-dynamical reflection equation (SDRE):
\be
A_{12}(\lambda) K_1(\lambda) B_{12}(\lambda) K_2(\lambda+\gamma h_1) = K_2(\lambda) C_{12}(\lambda)
K_1(\lambda+\gamma h_2) D_{12}(\lambda)\label{SDRE}
\ee
where $A, B, C, D$ are $c$-number matrices in $\mbox{End}(V) \otimes \mbox{End}(V) (\otimes {\mathbb C}[[u_1,\ u_2]])$ depending on the dynamical variables $\l=\{{\l}_i\}_{i\in\{1\cdots N\}}$ and possibly on spectral parameters, this last dependance being then encoded in the labeling $(1,2)$.
When one considers (as in \cite{AR})  non-operatorial or so-called ``scalar'' solutions
(i.e. dimension-1 representations
of the algebra) this $c$-number solution matrix will be denoted $k(\lambda)$.
The exact meaning of the shift on these dynamical variables $\l$ in~(\ref{SDRE}) together with
the main definitions and properties concerning the SDRA will be given in the next section and in appendix (A).

The characteristic feature of the SDRA is that the integrable quantum Hamiltonians, obtained by the
associated trace procedure from a monodromy matrix, exhibit an explicit dependance on the shift operators
$\exp{{\pd}_i}$,  ($ {\pd}_i={{\pd}\over{\pd{\l}_i}}$).
In the case of the previously constructed dynamical reflection algebra known as "dynamical boundary algebra"~\cite{KuMu} however, such a dependance also arises but may altogether vanish when the basic scalar reflection matrix $k(\l)$ used to build the monodromy matrix is diagonal~\cite{NA}.
In the case of Gervais-Neveu-Felder dynamical quantum group, an explicit dependance also occurs but the commutation of Hamiltonians requires to restrict the Hilbert space of quantum states to zero-weight states under the characteristic Cartan algebra defining the dynamical dependance~\cite{ABB}.
No such restriction occurs here which singles out the SDRA as the most useful algebraic framework to formulate spin-chain type systems with an extra potential interaction between the sites of the spins and explicit dynamics on the positions, on the line of the spin-Ruijsenaar-Schneider systems~\cite{KRI}.
Explicit formulae for these Hamiltonians, deduced in~\cite{NA} in the most generic frame, 
yield complicated--looking objects with intricated connections between spin interactions 
and "space-like" potential interactions.
Such formulations may however simplify, as shall be shown here,
when the building matrices $A, B, C, D$ take some particular form.

Our purpose here is twofold. In order to construct
consistent sets of $A, B, C, D$ structure matrices
we formulate generalized Yang-Baxter-type consistency equations (YBCE) 
extending the ones found in e.g.~\cite{FrMa,KuSkl}
with the same assumption of associativity of the SDRA.
This larger set of sufficient conditions is denoted 
``$g$-extended Yang-Baxter-type consistency equations''($g$-YBCE) 
since they depend upon an automorphism $g$ of the auxiliary space ${\cal V}$.
Analyzing and solving at least partially these two sets of equations (YBCE and $g$-YBCE) for the matrices $A, B, C, D$, we propose explicit parametrizations of the matrices $A, B, C, D$, the scalar solutions $k(\l)$
and the generating matrix $K(\l)$ in terms of quantum group-like algebraic structures ($R$ matrices
and Drinfel'd twists).

In a second step we plug these parametrizations into the general formulae for monodromy matrices, and obtain simplified expressions for them. These factorized expressions in terms of non--dynamical $R$ matrices and Drinfel'd twists,
simplify considerably the monodromy matrices found in \cite{NA} and represent therefore
a suitable starting point to construct
and solve quantum integrable Hamiltonians by allowing an explicit realization of the intricate formulae
previously obtained in \cite{NA}.
We also expect that this procedure may help to understand the nature of the algebraic structure implied by SDRE~(\ref{SDRE}), specifically its possible connections with ordinary quantum group structures through
Drinfel'd twists.
However we must emphasize that at every stage, including the all-important first step of deriving Yang-Baxter-type consistency equations, but also $A, B, C, D$  parametrizations, resolution of  (\ref{SDRE}) for non-operatorial $k(\l)$ matrix, and even comodule structure yielding the monodromy matrix, we have proceeded by sufficient conditions;
therefore we shall not cover here the full description of the 
algebraic content of~(\ref{SDRE}).

Our paper goes as follows.
In a first section we describe the notations, derive the 
sufficient Yang-Baxter-type consistency equations considered here, and discuss the 
possible factorization byof dynamical dependence in one of the four coefficient matrices.
A second section treats the case of the simplest set ($g = \1$)
of Yang-Baxter-type consistency equations, ending with the factorization of the monodromy matrix.
We develop in an already extensive way the 
analysis of this set of YBCE, in order to establish clearly in a first stage
the major steps of the
parametrization procedure, and the subsequent derivation of the factorized form
of the monodromy matrices, without the added complications induced
by the existence of a non--trivial automorphism. We also
discuss --more or less sketchily-- some alternative paths
to constructing different sets of solutions, by relaxing or eliminating
some of the restrictions defining our sufficient conditions.
A third section then deals with the full set of Yang-Baxter-type consistency equations
for a generic $g$. The main features remain, but the occurence of $g$ induces several subtle
effects, and requires the introduction of some supplementary
assumptions, which we discuss in detail.
Finally some conclusions and perspectives are drawn.

\section{Notations and derivation of the two sets of Yang-Baxter type consistency equations}

The main features of the reflection equations yielding the SDRA are given in appendix (A).
In this section we will thus start with the SDRE~(\ref{SDRE}), recall the definitions and properties of the objects it involves and obtain two sets of Yang-Baxter-type consistency equations (YBCE and g-YBCE).

We start by expliciting the exact meaning of the shift on the so-called dynamical variables $\l$ in~(\ref{SDRE}).
Let $\mf{g}$ be a simple complex Lie algebra and $\mf{h}$  a  commutative subalgebra of $\mathfrak{g}$ of dimension $n$.
(For an extension to non-commutative  $\mf{h}$ see \cite{Ping}.)

Let us choose a basis $\ds\{h^{i}\}_{i=1}^n$ of $\mf{h}^{\ast}$ and
let $\ds \l =\sum_{i=1}^n\l_i h^i$, with $(\l_i)_{i\in \{1,\cdots, n\}}\in \C^n$ be an element of $\mf{h}^{\ast}$.
The dual basis is denoted in $\mf{h}$ by $\ds\{h_{i}\}_{i=1}^n$ \footnote{This
labeling of the dual basis must not be confused with the traditional labeling
of auxiliary spaces in the global formulation of the SDRE~(\ref{SDRE}) }.
For any differentiable function
$f(\l)=f(\{\l_i\})$ one defines:

\begin{eqnarray}
f(\l+\g h)= e^{\g \mathcal{D}}f(\l)e^{-\g \mathcal{D}},
\;\mbox{where} \;\;
\mathcal{D}=\sum_i h_i \partial_{\l_i} \nn
\end{eqnarray}

It can be seen that this definition yields formally

\begin{eqnarray}
f(\lambda+\gamma h)=  f(\{\lambda_i+\gamma h_i\}) = \sum_{m \geq 0}\frac{\gamma ^m}{m!} \sum_{i_1, \dots ,i_m = 1}^n
\frac {\partial^m f(\lambda)}{\partial \lambda_{i_1} \dots
\partial \lambda_{i_m}}  h_{i_1} \dots h_{i_m} \nn
\end{eqnarray}

which is a function on $\mf{h}^{\ast}$ identified with $\C^n$ taking values in $\mathsf{U}(\mathfrak{h})$.

From now on, in order to alleviate the notations, we shall denote $f(h) \equiv f(\lambda+\gamma h)$.

Assumption of the associativity of SDRA and comparison of two possible ways of exchanging three $K$ matrices
requires zero weight conditions on structure matrices, namely:
\be
\left[h_i \otimes \1, B_{12} \right]=0,\quad
\left[\1 \otimes h_i ,C_{12} \right]= 0,\quad
\left[h_i \otimes \1+ \1 \otimes h_i ,D_{12} \right] =0,
\quad \forall i\in \{1,\cdots, n\}
\ee

It then yields the Yang--Baxter type consistency equations.
A derivation of such sufficient consistency conditions yielding the YBCE is found (for the non-dynamical case) in e.g. 
\cite{FrMa} and for the semi-dynamical  case in \cite{NAR}.

Here this derivation yields the following set of
Yang--Baxter equations

\begin{equation}\label{ybes}
\begin{array}{l}
a)\quad A_{12}\ A_{13}\ A_{23} =A_{23}\ A_{13}\ A_{12}
 \\
b)\quad A_{12}\ C_{13}\ C_{23} = C_{23}\ C_{13}\ A_{12}(h_{3})
\\ 
c)\quad D_{12}\ B_{13}\ B_{23}(h_{1})= B_{23}\ B_{13}(h_{2})\ D_{12}
\\
d)\quad D_{12}(h_{3})\ D_{13}\ D_{23}(h_{1})= D_{23}\ D_{13}(h_{2})\ D_{12}
\end{array}
\end{equation}

This set, obeyed for instance by the constant (i.e. non-spectral parameter dependant) $A, B, C, D$ matrices~\cite{ACF} associated to the  Ruijsenaar-Schneider (RS) $A_{n}$ rational and trigonometric models~\cite{RS}, will be globally denoted as ``standard Yang--Baxter type consistency  equations'' or YBCE. It is in fact the simplest example of a more generic form derived presently, but it is worth separating it in our derivation of parametrization of solutions, so as to treat it
as a first simpler example even though it already exhibits the essential features of this parametrization.

A more general form of Yang--Baxter type consistency equations is indeed derived from
(\ref{ybes}) once one notices that the
identification of the structure matrices $A, B, C, D$ in
(\ref{SDRE}) exhibits some freedom due to the invariance of (\ref{SDRE}) under suitable transformations. In particular, the exchange algebra encapsulating the exchange relations for the generators of the SDRA building the matrix $K$ (understood as an object in $\mbox{End}{\cal V} \otimes \mf{a}$ where $\mf{a}$ is the SDRA) can be equivalently formulated by multiplying the l.h.s. of (\ref{SDRE}) by $g \otimes {\1}$, where $g$ is an automorphism of the auxiliary space ${\cal V}$ (see Appendix A for notations on the auxiliary space). 

{\bf Remark:} The complete multiplication of~\ref{SDRE} by {\it two}
automorphisms $g \otimes g'$ can always be brought back to this form by a global change of basis on ${\cal V}$
parametrized by $g'$, multiplying the r.h.s. of
(\ref{SDRE}) by ${g'}^{-1}\otimes {g'}^{-1}$, for ${g'}$ any automorphism on ${\cal V}$. The endomorphisms $h$ representing the generators of the Lie algebra ${\mf h}$ acting on ${\cal V}$ (assumed to be a diagonalizable module of ${\mf h}$) are accordingly redefined as $g'hg'^{-1}$.

In order to be able to undertake some specific technical manipulations,
we shall restrict $g$, in the case when ${\cal V}$ is
an evaluation module with spectral parameter $u$, by requiring that
its adjoint action on any matrix in ($\mbox{End}V^{\otimes N} \otimes {\mathbb C}[[u_1...u_N]]$) yields again a ``factorized'' matrix in ($\mbox{End}V^{\otimes N} \otimes {\mathbb C}[[u_1...u_N]]$).
In other words the adjoint action of $g$ must be compatible with the evaluation representation.
This is equivalent to asking that, provided that $g$ admits an operatorial logarithm $\g=\mbox{log}g$,
$[[\g,u.],u.]=0$ where $u.$ is the automorphism of formal multiplication by $u$ on ${\cal V}$.
As an example, any automorphism $\g$ commuting directly with $u$ will provide a suitable $g=\mbox{exp}\g$.

We shall also be later interested in particularizing endomorphisms $\g$ such that $[\g,u.]=0$.
This is indeed equivalent to assuming that the action of $\g$ on ${\cal V}=V \otimes \C[[u]]$
is represented by afunctional matrix $M(\g)\in\mbox{End}V \otimes \C[[u]]$ acting on ${\cal V}$.
Such endomorphisms will be called "factorizable" for obvious reasons.
Automorphisms of the type $g=\mbox{exp}\g$ with $[[\g,u.],u.]=0$ will be called "ad-factorizable".

This l.h.s. gauging of (\ref{SDRE}) now leads to a new definition of structure matrices: \be \tilde A_{12} = g_{1}\ A_{12}\ g_{2}^{-1},~~~~~\tilde B_{12} = g_{2}\ B_{12}, ~~~~~\tilde C_{12}= g_{1}\ C_{12}, ~~~~~\tilde D_{12} = D_{12}. \label{gauge} \ee If we now assume perfectly consistently that $\tilde A, \tilde B, \tilde C, \tilde D$ (instead of $A, B, C ,D$) obey the sufficient equations (\ref{ybes})) we get a new set of Yang--Baxter type
consistency equations for $A, B, C, D$:
\begin{equation}\label{yben}
\begin{array}{l}
a)\quad A_{12}\ A_{13}^{gg}\ A_{23} = A_{23}^{gg}\ A_{13}\ A_{12}^{gg} 
\\ 
b)\quad A_{12}\ C_{13}^{g_{1}}\ C_{23} = C_{23}^{g_{2}}\ C_{13}\ A_{12}^{gg}(h_{3})
\\ 
c)\quad D_{12}\ B_{13}\ B_{23}^{g_{3}}(h_{1})= B_{23}\ B_{13}^{g_{3}}(h_{2})\ D_{12}
\\
d)\quad D_{12}(h_{3})\ D_{13}\ D_{23}(h_{1})= D_{23}\ D_{13}(h_{2})\ D_{12}
\end{array}
\end{equation}
where $X_{12}^{g1}$, $X_{12}^{g2}$ and $X_{12}^{gg}$ now denote respectively the following adjoint actions 
$g_{1} \ X_{12}\ g^{-1}_{1}$,
$g_{2}\ X_{12}\ g_{2}^{-1}$ and
$g_{1}\ g_{2}\ X_{12}\ g^{-1}_{1}\ g_{2}^{-1}$.

The generating matrix $K$ is unmodified under this operation, and will thus be used directly when building monodromy matrices from the comodule structure.
Consistency however will require to use tilded matrices
(\ref{gauge}) to build the $N$-site monodromy matrix. This set of equations is hereafter denoted ``$g$-deformed Yang--Baxter type consistency equations'' or $g$-YBCE.

It is interesting to note that although the tilded ``structure matrices'' are not obtained as adjoint actions of $g$ on the $c$-number original matrices $A, B, C$, and may therefore not be represented as finite-size
matrices in the evaluation representation when ${\cal V} = V \otimes \C[[u]]$,  the new Yang--Baxter equations exhibit only adjoint actions of the automorphism $g$ on the original $c$-number matrices $A,B,C$, hence are again written in terms of finite-size numerical
matrices as follows from our restriction on $g$. On the example in \cite{ACF} where $g = \exp[{d \over d u}]$, $u$ being the spectral parameter in the evaluation representation on ${\cal V} = V \otimes \C[[u]]$, it appears that in this case, although the structure matrices (\ref{gauge}) are not $c$-number matrices anymore (in other words, $V \otimes \C[[u]]$ is not an evaluation
module for (\ref{gauge})) the Yang--Baxter equations themselves admit a representation (\ref{yben}) on the evaluation module, allowing the normal matrix manipulations to parametrize its solutions. Auxiliary action is here a shift of the spectral parameter.

We shall impose two further restrictions on $g$. The first is purely
technical: we shall assume the existence of an endomorphism $\log\ g$ on ${\cal V}$ such that $\exp [\log\ g] = g$. This will be used later when solving the so-called quasi-non dynamical conditions on given matrices acting on ${\cal V}$ or ${\cal V} \otimes {\cal V}$. The second one will impose that $g$ does not depend on dynamical variables; it will play a central role when solving the Yang--Baxter equations.

It is finally relevant to start at once discussing the possible
parametrizations of the $D$ matrix which can essentially be
treated (as will be seen in the next sections) independently of $A, B, C$. 
Analyzing the possibilities of existence of invertible scalar (non-operatorial) solutions $k(\lambda)$ to (\ref{SDRE})
leads us to consider three possible situations for the relevant parametrizations of $D$. They will take a general form:
\be D_{12} = q_{1}^{-1}\ q_{2}^{-1}(\lambda +h_{1})\ \tilde R_{12}\ q_{12}(\lambda+h_2)\ q_2 \label{paramt} \ee where $q$ is a scalar dynamical matrix in $\mbox{End} V$ or $(\mbox{End} V)\otimes \C[[u]]$ (factorizable). The three possibilities to consider are the following:
\\
\\
1. Existence of decomposition (\ref{paramt}) with a non-dynamical $R$-matrix $\tilde R$:
\be \tilde R_{12}(\lambda+h_3)=\tilde  R_{12}(\lambda) \nonumber \\ \tilde  R_{12}\ \tilde R_{13}\ \tilde R_{23} = \tilde R_{23}\ \tilde R_{13}\ \tilde R_{12} \label{paramt1} \ee
\\
2. Existence of a decomposition (\ref{paramt}) with a quasi-non dynamical $R$-matrix i.e. \be \tilde R_{12}
(\lambda+h_3)= f_{1}\ f_{2}\ \tilde R_{12}(\lambda)\ (f_{1})^{-1}\ (f_{2})^{-1} \nonumber\\ R_{12}^0\ R_{13}^{0 f_{1} f_{3}}\ R_{23}^0= R_{23}^{0 f_{2} f_3}\ R_{13}^0\ R_{12}^{0 f_1 f_2} \nonumber\\ \mbox{where} ~~~~~R_{12}^0 = Ad. exp[-\sigma(\log f_1 +\log f_2)]\ \tilde R_{12}(\lambda) \nonumber\\ \mbox{so that} ~~~~~R_{12}^0(\lambda+h_3) = R_{12}^0 ~~~~\mbox{non-dynamical}.  \label{paramt2}\ee here $f$ is an ad-factorizable automorphism of ${\cal V}$, not necessarily identified with the automorphism $g$ in (\ref{yben}).
\\
3. Neither decomposition exists.
\\
\\
{\bf Remark}: Situations 1 and 2 may coexist, but we shall not establish if and when such a coexistence arises, it being not relevant for our specific purpose.

Possibility 1 (hereafter denoted ``de-twisting of the $D$ matrix'') is indeed realized when the $D$-matrix is the representation of the universal $R$ matrix for the quasi--Hopf algebra obtained by Drinfel'd twist of a Hopf algebra. $\tilde R$ is then the representation of the universal $R$-matrix for the Hopf algebra \cite{Drin,JKOS,ABRR}.
By extension of this notion we shall sometimes denote as ``twist'' the shifted conjugation by $q$
in (\ref{paramt}) and ``twist matrix'' the $q$ matrix.

It was recently proven \cite{BRT} at the level of universal $R$ matrices that $D$- matrices of weak
Hecke type, associated to the $A_{n}$ simple Lie algebra, could always be constructed as Drinfel'd twists of non-dynamical Cremmer--Gervais \cite{CG} $R$-matrices \be D_{12} =g_{1}^{-1}\ g_{2}^{-1}(h_{1})\ R_{12}^{CG}\
g_{1}(h_{2})\ g_{2}. \label{drin}\ee However, even in the case  of simple $A_n$ Lie algebra (no spectral parameter) exhaustive resolutions of the dynamical Yang--Baxter equation shows that non-weak-Hecke type solutions exist \cite{AR2}. In addition, the case of $A_{n}^{(\ref{SDRE})}$ affine Lie algebra (naturally relevant when $D$ depends on a spectral parameter) is not covered by the result in \cite{BRT}. We shall hereafter be lead to differentiate between the cases where $D$ can be ``detwisted'' as in (\ref{drin}), and cases
where $D$ can not be written as in (\ref{drin}). This is in
particular relevant to study the possible existence and precise constructions of invertible $c$-number solutions $k(\lambda)$.

Possibility 2 (hereafter denoted ``quasi--detwisting of $D$-matrix) has as far as we know no such interpretation yet, but should have a relation with the Drinfeld twist formulation in the context of the $g$--deformed Yang--Baxter equations.

We can now start the discussion on parametrization of $A, B, C, D$ and $K$ and construction of monodromy matrices and Hamiltonians, starting with the simpler case of standard Yang--Baxter type equations (\ref{ybes}).

\section{Standard Yang-Baxter type consistency equations}

\subsection{The $A,\ B,\ C$ matrices}

Once again ${\cal V}$ is either a finite dimensional vector space $V$ or an evaluation module $V \otimes {\mathbb C}[[u]]$.
We assume that the vector space $V$ is an irreducible representation of the dynamical algebra ${\mathfrak h}$.
Since $B_{12}$ is a space-1 zero weight matrix, and choosing from now on ${\mf h}$ to be the
Cartan algebra of $(gl(n))$, $B$ can be parametrized as
\be B = \sum_{i=1}^{n} e_{ii}\otimes b_{i}(\lambda)~~~~~b_{i}(\lambda) \in \mbox{End} V \otimes {\mathbb C}[[u]], \ee
Since $D$ is a zero-weight matrix, it can be parametrized as
\be D = \sum_{i,j =1}^n d_{ij}(\lambda) \; e_{ii} \otimes e_{jj} + \sum_{i \neq j =1}^n \Delta_{ij}\; (\lambda) e_{ij} \otimes e_{ji}. \label{D} \ee

Equation (\ref{ybes}c) now reduces to \be d_{ij}\ (b_{i}\ b_{j}(h_{i})- b_{j}\ b_{i}(h_{i})) =0. \label{constr} \ee
We shall from now on, until the end of the paper, assume
that all diagonal elements $d_{ij} \neq 0$, for all $i,j \in\{1, \ldots,n \}$.

In this case $b_{i}\ b_{j}(h_{i}) = b_{j}\ b_{i}(h_{i})$ for all $i,\ j$. If all $b_i$'s are invertible ($n\times n$) matrices, this implies that $b_{i}$ are parametrized as: 
\be b_{i} = b^{-1}\ b\ (\lambda_{j} +\delta_{ij} \gamma) \quad \mbox{with} \; b \; \mbox{some invertible matrix} \label{param} \ee 
or equivalently
\be B_{12} = {\1} \otimes b^{-1}\ b(h_{1})=b_{2}^{-1}\ b_{2}(h_{1})\label{param2} \ee 
using the compact dynamical shift notation and space indices. Here again 
$b(\lambda) \in \mbox{End} V \otimes {\mathbb C}[[u]]$.

If some $b_i$'s are not invertible the simple parametrization
(\ref{param}) is not available. Examples of such situations are easily given. Define a set of $n$ mutually commuting projectors $P_{i}$, such that in addition $[P_{i},\ b] =0$, then \be b_i = P_i\ b^{-1}\ b(\lambda_j +\delta_{ij}\gamma) \label{param3} \ee obeys (\ref{constr}). It is not clear however whether an exhaustive classification of such solutions may be available.

If $B$ is invertible, plugging back $C = B^{\pi}$ into
(\ref{ybes}b) yields the simple identity: \be  (b_{1}\ b_{2}\ A_{12}\ b_{1}^{-1}\ b_{2}^{-1})(h_{3}) =  b_{1}\ b_{2}\ A_{12}\ b_{1}^{-1}\ b_{2}^{-1} \label{ident} \ee equivalently stating that $b_{1}\ b_{2}\ A_{12}\ b_{1}^{-1}\ b_{2}^{-1} =R_{12}$ is non-dynamical. Furthermore plugging it into (\ref{ybes}a) immediately implies that $R_{12}$ is a non-dynamical solution of the Yang--Baxter equation, or a non-dynamical $R$ matrix.

If $B$ is non-invertible, the absence of explicit parametrization prevents us from deriving a general form for $A$. However the example (\ref{param3}) for instance is workable. Defining once again: $~R_{12} = b_{1}\ b_{2}\ A_{12}\ b_{1}^{-1}\ b_{2}^{-1}$ yields from (\ref{ybes}b) 
\be R_{12}\ P_{i1}\ P_{i2} = P_{i1}\ P_{i2}\ R_{12}(h_{i}) \label{ybenp} \ee 
and from (\ref{ybes}a) again the Yang--Baxter equation for $R$. Once again it may not be possible
to exhaust all simultaneous solutions to Yang--Baxter equations and (\ref{ybenp}). However one deduces that if $R$ is a non-dynamical $R$ matrix and $\{ P_{i}\}$ a set of projectors such that $[P_i \otimes P_i ,\ R] =0$ and $[P_{i},\ b] =0$ then they provide a consistent set of matrices  \be && A_{12} = b_{1}\ b_2\ R_{12}\ b_{1}^{-1}\ b_{2}^{-1} \nn\\
&& B =C^{\pi} = \sum e_{ii} \otimes P_{i}\ b^{-1}\
b(\lambda_{i}+\gamma). \label{gen}\ee Such projectors exist e.g. if $R$ is a Yangian-type solution in $A_{n}^{(\ref{SDRE})} \otimes A_{n}^{(\ref{SDRE})}$ \be R  = {\1} \otimes {\1} + {\Pi_{12} \over \lambda -\mu} \label{yang} \ee since then for {\it any} projector $[P\otimes P,\ R] =0$. Choosing these projectors $P$ to commute with an arbitrary chosen matrix $b$, and with each other (e.g. elements among the set of projectors on eigenvectors of $b$) one gets $A, B,$ and $ C$.

To conclude: If $d_{ij} \neq 0$ for all $i,\ j$, and $B$
invertible there exists a parametrization of $A, B, C$ as: \be && A = b_{1}^{-1}\ b_{2}^{-1}\ R\ b_{1}\ b_{2} \nn\\ && B = C^{\pi} = {\1} \otimes b^{-1}\ b(\lambda+h_{1}) = b_{2}^{-1}\ b_{2}(\lambda+h_{1})\label{param01}\ee where $R$ is a non-dynamical quantum $R$-matrix and $b$ some dynamical matrix.\\
One immediately establishes here:
\\
\\
{\bf Proposition 2}
\\
If $A, B, C$ are parametrized as in (\ref{param01}) by matrices $b$ and $R$, the following statements are equivalent

\begin{itemize}

\item{(a) The SDYBE equation (\ref{SDRE}) has an invertible scalar solution $k(\lambda)$}

\item{(b) $D$  can be de-twisted, following (\ref{drin}), to a non-dynamical matrix $R$ with twist given by $q = bk$.}

\end{itemize}

{\bf Proof:}
\begin{itemize}

\item{(a) $\Rightarrow$ (b) by direct inversion of (\ref{SDRE}) yielding (\ref{drin}) with $q = bk$}

\item{(b) $\Rightarrow$ (a) by direct plug of (\ref{drin}) into (\ref{SDRE}) using $q$ and $b$ yielding $k = b^{-1} q$ as a scalar solution.}

\end{itemize}
Hence, whether $D$ can not be detwisted at all or can not be
detwisted to $R$ the absence of a scalar invertible solution may cause serious practical issues to build integrable spin-chain type Hamiltonians. However, if $D$ is de-twistable to another $\tilde R$, one may nevertheless draw interesting conclusions regarding possible non-invertible scalar solutions, and even monodromy matrices. We shall henceforth proceed with our general trichotomy.

\subsection{The $D$ matrix and $K$ solutions}

As indicated above, we shall separate this discussion into three subcases, whether or not $D$ can be detwisted as in (\ref{drin}) and whether it is detwisted as in (\ref{paramt1}) or (\ref{paramt2}). Note immediately that one can show easily:

\subsubsection{Case 1 and 2. D is detwistable or quasi-detwistable}

We use here the general form

\be D_{12} = q_{1}^{-1}(\lambda +h_{2})\ q_{2}^{-1} \tilde R_{12}\
q_{1}\ q_{2}(\lambda +h_{1}) \label{param02} \ee  where $\tilde R$ is either non-dynamical or quasi non-dynamical. If $A, B, C$ are parametrized as in (\ref{param01}), plugging (\ref{param01}) and (\ref{param02}) into (1.1) leads to the following equation \be R_{12}\ (bKq^{-1})_1\ q_{1}\ (bKq^{-1})_{2}(h_{1})\ q_{1}^{-1} = (bKq^{-1})_{2}\ q_{2}\ (bKq^{-1})_{1}(h_{2})\ q^{-1}_{2}\ \tilde R_{12}. \label{geq} \ee General solutions to (\ref{geq}) are not obvious to formulate due to the coupling between spaces 1 and 2 induced by the adjoint action of $q_{1,2}$ on $(b\ k\ g^{-1})_{2,1}(h_{1,2})$. If however $b\ K\ q^{-1}$ is such that: 
\be(b\ K\ q^{-1})_{1}(h_{2}) = {\cal A} \otimes {\1}={\cal A}_{1}
\label{cond} \ee 
for some matrix functional ${\cal A}$  then (\ref{geq}) simplifies to a Yang--Baxter type form \be R_{12}\ (b\ K\ q^{-1})_{1}\ {\cal A} (b\ K\ q^{-1})_{2} = (b\ K\ q^{-1})_{2}\ {\cal A} (b\ K\ q^{-1})_{1}\ \tilde R_{12}. \label{tybe} \ee

Condition (\ref{cond}) can be explicitly solved as follows. From the general definition of shifts, applied to the $gl(n)$ case, one has \be (b\ K\ q^{-1})_{1}(h_{2}) = \sum _{i=1}^{n} b\ K\ q^{-1}(\lambda_{j} +\gamma \delta_{ij}) \otimes e_{ii}. \label{cond2} \ee Factorizing ${\1}_{2}$ as in (\ref{cond}) requires to have \be b\ K\ q^{-1}(\lambda_{j} +\gamma \delta_{ij})= b\ K\ q^{-1}(\lambda_{j} +\gamma \delta_{lj}) \label{cond3} \ee for any index pair $(i,\ l)$. This is equivalent to restricting $b\ K\ q^{-1}$ to depend on the following new set of dynamical variables \be \sigma = \sum_{i=1}^{n} \lambda_{i}, ~~~~~\theta_{i} = \sigma -2\lambda_{i}, ~~~~~i =2, \ldots, n \label{new} \ee constrained by : $b\ K\ q^{-1}(\theta_{i} +2\gamma) = b\ K\ q^{-1}(\theta_i)$ for $i = 2, \ldots, n$.

Equation (\ref{geq}) now becomes a usual dynamical Yang--Baxter intertwining equation 
for $\kappa \equiv bKq^{-1}$for the simplified situation where $R$ itself is non-dynamical \be R_{12}\ \kappa_{1}(\sigma)\ \kappa_{2}(\sigma +\gamma) = \kappa_{2}(\sigma)\ \kappa_{1}(\sigma +\gamma)\ \tilde R_{12}. \label{nondy} \ee We shall not discuss (\ref{nondy}) in full generality. We now separate our discussion into
two subcases. 

\subsubsection{D detwistable, $\tilde R$ non--dynamical}

Two simple and relevant examples will now provide us with explicit realizations of solutions to the SDRE (\ref{geq}).
\\
\\
{\bf a) Non dynamical quantum group}

Given any non-dynamical solution ${\cal Q}$ to: \be R_{12}\ {\cal Q}_{1}\ {\cal Q}_{2} = {\cal Q}_{2}\ {\cal Q}_{1}\ \tilde R_{12} \label{nondy1} \ee $K(\lambda) = b^{-1}\ {\cal Q}\ q(\lambda)$ realizes a solution of (\ref{geq})). 
In particular if ${\cal Q}$ is a factorized matrix, represented in $\mbox{End} V \otimes {\mathbb C}[[u]]$,
$K(\lambda)$ is also such a solution to (1.1).
It follows that:
\\
\\
{\bf a1.} if $R =\tilde R$ ($\leftrightarrow$ existence of scalar
invertible solution)
\\
\\
Any realization ${\cal Q}$ of the quantum group described by the RTT
formulations with $R$ as evaluated $R$ matrix, will provide a realization of the SDRA as $K= b^{-1}\ {\cal Q}\ q$. This includes as well scalar solutions (yielding scalar $k$ matrices) or operator like solutions (representations of the quantum group by operators on some Hilbert space ${\cal H}$). In particular, ${\cal Q}=\1$ yields an invertible scalar solution $k = b^{-1}\ q$, consistent with Proposition 2.
\\
\\
{\bf a2.} if $R \neq \tilde R$ (no invertible scalar solutions)
\\
\\
Then any intertwiner matrix (scalar or operational) ${\cal Q}$: \be
R_{12}\ {\cal Q}_{1}\ {\cal Q}_{2} = {\cal Q}_{2}\ {\cal Q}_{1}\ \tilde R_{12} \label{nondy2} \ee provides us with realizations of the SDRA.
\\
\\

{\bf b) Quasi-non dynamical quantum group}

Let us consider the more general quadratic exchange relation: 
\be R_{12}\ {\cal Q}_{1}\ (a\ {\cal Q} a^{-1})_{2} = {\cal Q}_{2}\ (a\ {\cal Q} a^{-1})_{1}\ \tilde R_{12} \label{mgen} \ee 

for some ad-factorizable automorphism $a$ of the auxiliary space 
${\cal V}$, such that $[a \otimes a,\ R] = [a\otimes a,\ \tilde R] =0$. From any non-dynamical representation 
${\cal Q}$ of this exchange algebra (scalar or operatorial) one can build a representation (scalar or operatorial) of the SDRA as: \be K =
b^{-1}(\lambda)\ (\exp[\sigma \log a]\ {\cal Q}
\exp[-\sigma \log a]) q(\lambda). \label{dysol} \ee assuming the existence of a logarithm of $a$. This adjoint action transforms the dynamical shift on any dynamical parameter $\lambda$ into a conjugation by $a$, yielding what we 
will call quasi-non dynamical condition for $\tilde q = (\exp[\sigma \log a]\ {\cal Q}
\exp[-\sigma \log a])$. 

\be \tilde q (\lambda +h_{2}) = a\ \tilde q(\lambda)\ a^{-1} \otimes {\1}. \label{quasi} \ee

Once again, ad-factorizability of $a$ ensures that (\ref{mgen}) and (\ref{dysol})
are finite-matrix algebraic equations on the auxiliary space $V$.

\subsubsection{D quasi-detwistable, $\tilde R$ quasi-non dynamical}

Here one assumes that $\tilde R$ obeys (\ref{paramt2}) for some ad-factorizable automorphism $f$
of ${\cal V}$. It is still possible to obtain explicit representations of (1.1) 
as modified versions of the representations  given in the previous subsection. Namely the 
non--dynamical quantum group (NDQG) construction {\bf a)} is modified as follows: (\ref{nondy1})
becomes \be R_{12}\ {\cal Q}_{1}\ {\cal Q}_{2}f_2= {\cal Q}_{2}\ {\cal Q}_1f_1\ \tilde R_{12}^0 \label{mod1}  \ee where $\tilde R^0$ is the non-dynamical part of $\tilde R$ 
extracted from (\ref{paramt})
\be \tilde R_{12}(\lambda) = Ad. exp[-\sigma (\log f_1 + \log f_2)] \tilde R^{0}_{12}\label{mod2} \ee and $K(\lambda)$ becomes
\be K(\lambda) = b^{-1}(\lambda)\ {\cal Q}\ exp [-\sigma \log f]\ q(\lambda). \label{mod3} \ee The Quasi-NDQG
{\bf b)} is modified as follows: (\ref{mgen}) becomes \be R_{12}\ {\cal Q}_{1}\ (a
{\cal Q} a^{-1})f_{2}\ = {\cal Q}_{2}\ (a {\cal Q} a^{-1})f_{1}\ \ \tilde R_{12}^0 \label{mod3b} \ee
with the $K$ matrix now being \be K = b^{-1}(\lambda) \Big(Ad. exp[\sigma \log a]\ {\cal Q} \Big )exp[\sigma \log f]\ q(\lambda). \label{mod4} \ee Note that here no relation between the two automorphisms
$a$ and $s$ need be assumed. However if $f$, although ad-factorizable, is not factorizable (see e.g. $\log f \equiv  d/du$),  equations (\ref{mod1}) and (\ref{mod3b}) cannot be written as algebraic equations for finite-size
matrices in $\mbox{End} V \otimes {\mathbb C}[[u]]$, and the objects ${\cal Q}$, solutions of
(\ref{mod1}) and (\ref{mod3b}), may not be expandable
in formal power series of the variable $u$; subsequent interpretation of $K$ as a generating functional for
some quantum algebra is then unavailable, and the correct interpretation of (\ref{SDRE}) in this
context remains to be explicited.

\subsubsection{Case 3. $D$ non de-twistable}

One is here able to build new sets of realizations $K(\lambda)$ of  the SDRA if one knows at least one (non invertible!) scalar solution $K(\lambda)$, from the left comodule structure, described as follows:
\\
\\
{\bf Proposition 3}
\\
If $K_{0}(\lambda)$ is a solution of (\ref{SDRE}), and $A, B, C$ are
parametrized by $R$ and $b$ as in (\ref{param01}), from any
solution of: \be R_{12}\ q_{1}\ q_{2}(h_{1})  = q_{2}\
q_{1}(h_{2})\ R_{12} \label{qq} \ee 
such that, once again,
$q_{n}(h_{m})= {\cal A}(q)_{n}\otimes {\1}_{m}$ (indices $n$ and
$m$ refer here to the labeling of auxiliary spaces in a multiple tensor
product),  one can build a
solution $b^{-1}\ q\ b\ K_{0}$ of (\ref{SDRE}).

One recovers once again the equations in (\ref{nondy1}) or (\ref{mgen}) (for $R=\tilde R $). Any (scalar) solution to (1.1) can be dressed to another solution, using any representation of the quantum group, or even quasi non-dynamical 
quantum group.

However when $D$ can not be detwisted, one cannot simplify the formulation of the monodromy matrix derived even from the simplest comodule structure of SDRA, hence we shall not consider this case in the next section.

\subsection{Monodromy matrices}

We shall restrict ourselves to the case where $A, B, C$ are
parametrized by matrices $R$ and $g$ (no $d_{ij} =0$), and $D$ is detwistable to a non--dynamical
$R$ matrix. In addition we shall only construct the monodromy matrix corresponding to
the simplest comodule realizations of the SDRA, i.e. realizations by $A, B, C, D$ matrices themselves (the
specific construction of new comodule realizations using the parametrizations derived here goes beyond the intended scope of this study). Moreover we shall also consider the simplest, i.e. non-dynamical, realizations of scalar $k$ matrices (\ref{nondy1}), (\ref{nondy2}). Construction of monodromy matrices to yield commuting spin-chain type Hamiltonians is mostly relevant from a physical point of view when the scalar solutions are themselves invertible. We shall nevertheless also consider the non-invertible, detwistable case as well, but once again only where $D$ is detwisted to a 
non-dynamical $R$ matrix.

\subsubsection{Existence of invertible solutions $k$}

We shall recall that one can then parametrize $A, B, C, D$ as 
\begin{equation} 
\begin{array}{l}\label{para} 
a) \quad B= C^{\pi}= {\1} \otimes b^{-1}\ b(\lambda +h_{1})
\\ 
b) \quad A= b_{1}^{-1}\ b_{2}^{-1}\ R \ b_{1}\ b_{2} 
\\ 
c) \quad D= b_{1} k_{1}(h_{2})^{-1}\ (b_{2} k_{2})^{-1}\ R \ b_{1} k_{1}\ b_{2} k_{2} 
\end{array}
\end{equation}
where $k$ is a particular invertible solution of (\ref{SDRE}). 
Other invertible solutions are given by: 
\be 
\tilde k = b^{-1}\ {\cal Q}\ b\ k,
\quad \mbox{where}\; {\cal Q} \;\mbox{is a scalar solution to} \quad R_{12}\ {\cal Q}_{1}\  {\cal Q}_{2} =  {\cal Q}_{2}\   {\cal Q}_{1}\ R_{12}. \label{inv2} 
\ee

There may be other invertible solutions obtained by a general resolution of (\ref{geq}), but at this stage we have no explicit parametrization for them and we shall therefore restrict ourselves to the previous dressed solutions $~b^{-1}{\cal Q}\ b\ k$.

We are now in a position to reformulate the monodromy matrix for a spin-chain type model, obtained from the particular comodule structure of the SDRA and the quantum trace structures detailed in \cite{NAR,NADR}, by plugging (\ref{para}),(\ref{inv2}) into the general formula. Denoting in addition by $\chi_{0}$ the solution to the dual 
SDRE required to build a ``reflection'' monodromy matrix, we recall that the $N$--site monodromy matrix can be chosen
of either two forms, to yield local Hamiltonians \cite{NA} by a
(partial) trace procedure over the finite vector space $V$ whichever structure is chosen
for the auxiliary space ${\cal V}$ \be
\chi_{0}^t\ A_{0\ 2N}\ C_{0\ 2N-1}\ldots A_{02}({h_<^{odd}})\
C_{01}\ T_{0}(h_<^{odd})\ D_{01}\ B_{02} \ldots D_{0\ 2N-1}\ B_{0\ 2N}\ e^{ {\cal D}_{0}} \label{mono} \ee or ($A \to C,\ B \to D$) making use of the  first known comodule structure. 
\\
Remark: the notation $X_{0a} (h_<^{odd})$ was introduced in \cite{NAR} and denotes 
$X_{0a} (\lambda + \Sigma_{n=0}^{ E(a/2)-1} h_{2n+1})$.

One may also use as ``site'' matrices $A \to (A^{-1})^{T},\ B \to (B^{-1})^T,\ C \to (C^{-1})^T,\ D \to (D^{-1})^T$ but we shall not consider this alternative possibility here for the sake of simplicity. Note
also the crucial occurrence of the shift operator $\exp[ {\cal D}_{0}]$ in the formulation of the monodromy ``matrix''. This guarantees that partial traces of monodromy matrices over the finite vector space $V$ commute as operators acting on the tensor product of the spin chain Hilbert space (in this case $({\mathbb C}^{n})^{\otimes N}$) and the functional space of differential functions over ${\mathfrak h}^*$. The price to
pay is that these traces lie not in the quantum reflection algebra
defined by (1.1), but in the extended operator space containing in addition derivatives w.r.t. variables in $({\mathfrak h}^*)^*$, such as built e.g. in \cite{Konno}. It may be conjectured that the relevant traces operate not in a quantum group but in a quantum groupoid structure relevant to the dynamical Yang--Baxter algebras \cite{Ping}.

The monodromy matrix (\ref{mono}) then becomes \be {\cal
O}_{N}^{-1}(\sigma) \Big \{ \chi^t_{0}\ b_{0}^{-1} R_{0\ 2N}\ \ldots R_{02}\  {\cal Q}_{0}\ R_{01} \ldots R_{0\ 2N-1}\ b_{0}\ k_{0}\ e^{\partial_{0}}\Big \}{\cal O}_{N}(\sigma) \label{mono1} \ee where the operator ${\cal O}_{n}(\sigma)$ acts only on the quantum spaces: \be {\cal O}_{N}(\sigma) = b_{2N}\ b_{2N-1}\ k_{2N-1}\ (b_{2N-2})(h_{2N-2}) \ldots b_{1} k_{1}(h_{3} + \ldots h_{2N-1}). \label{op}\ee

\subsubsection{No invertible solutions, $D$ detwistable to non-dynamical $\tilde R$ }

This corresponds to a situation where eq. (\ref{para}c) is replaced by \be D_{12}= q_{1}^{-1}(h_{2})\ q_{2}^{-1}\ \bar R_{12}\ q_{1}\ q_{2}(h_{1}) \label{p4}\ee but now $\bar R$ is a non-dynamical $R$-matrix not similar to $R$. In this case there exists no invertible scalar solution, otherwise $D$ could be detwisted to $R$. This situation is not so interesting from the point of view of realistic physical model building of spin chains, but it yields once again an interesting reduction of the monodromy matrix and may help in disentangling the general
structure of the semi-dynamical equation. Choosing the parametrization (\ref{para}a), (\ref{para}b), (\ref{p4}) and the scalar reflection solutions $\chi_{0}$ and $\tilde \chi_{0}$ one gets a
monodromy matrix: \be && {\cal O}^{-1}_{N}(\sigma)\Big \{\tilde
\chi_{0}^t\ b_{0}^{-1}\ R_{0\ 2N} \ldots R_{02}\Big (\prod_{k}^{1\to N}q_{2k+1}(h_{>}^{odd})\ b_{0}\chi_{0}(h_{>}^{odd})q_{0}^{-1}\ (\prod_{k}^{1\to N}q_{2k+1}(h_{>}))^{-1} \Big)\  \nn\\ && \bar R_{01} \ldots
\bar R_{0\ 2N-1}\ b_{0}k_{0}\ e^{\gamma \partial_{0}}\Big \}\ {\cal O}_{N}(\sigma) \label{mono2}\ee where \be {\cal O}_{N}(\sigma) = \prod_{k}^{1\to N}q_{2k-1}(h_{>}^{odd})\ b_{2k}(h_{>}^{odd}) \label{op2}\ee
 
 If $b_{0}\ \chi_{0}\ q_{0}^{-1}$ is non dynamical (i.e. if one chooses a solution $\chi_{0}$ of type
 given in subset 3.2.2a, a factorized compact formula for the monodromy matrix is then yielded with a form analogous to (\ref{mono1}).
However one must be careful that since no invertible scalar
solution $\chi_{0}$ to (\ref{SDRE}) exists, one has a priori no relation
expressing  a given dual solution $\tilde \chi_{0}$ in term of some direct solution $\chi$. 

This concludes our analysis of the semi-dynamical Yang--Baxter equation with ordinary Yang--Baxter conditions on $A, B, C, D$.

\section{g-deformed Yang--Baxter type consistency equations}

We shall for this discussion restrict ourselves to the simpler situation where all diagonal terms $d_{ij}$ of $D$ are non zero (as in Section 3), but also where matrices $B$ and $C$ are immediately assumed to be invertible. 
Once again, in the case where
${\cal V}$ is chosen to be an evaluation module ($\mbox{End}V\otimes {\mathbb C}[[u]]$)
we assume that the adjoint action
of the characteristic automorphism $g$ on any operator represented by 
a finite--size matrix in ($\mbox{End}V^{\otimes N} \otimes {\mathbb C}[[u_1...u_N]]$)
gives again a finite--size matrix (ad-factorizability condition).

\subsection{Parametrization of $A, B, C$}

Consider (\ref{yben}c) with the conditions $d_{ij} \neq 0$, B invertible. Eq. (\ref{constr}) is turned into \be b_{i}\ g\ b_{j}(\lambda_{i} +\gamma)\ g^{-1} = b_{j}\ g\ b_{i}(\lambda_{j} +\gamma)\ g^{-1} \label{constr2} \ee Assuming all $b_i$'s to be invertible , (\ref{constr2}) is solved by \be B_{12} = b_{2}^{-1}\ g_{2}\ b_{2}(\lambda+h_{1})\ g_{2}^{-1} \nn\\ C_{12}= b_{1}^{-1}\ g_{1}\ b_{1}(\lambda+h_{2})\ g_{1}^{-1}. \label{solu1}\ee

Define now \be R_{12} = b_{1}\ g_{2}\ b_{2}\ g_{2}^{-1}\ A_{12}\ (g_{1}b_{1}g_{1}^{-1})^{-1}\ b_{2}^{-1} \label{unt} \ee equation (\ref{yben}b) yields \be R_{12}(\lambda) = g_{1}\ g_{2}\ R_{12}(\lambda +h_{3})\ g_{1}^{-1}\ g_{2}^{-1} \label{solu2}\ee meaning that for any index $i$ \be R_{12}(\lambda_i + 1) = g_{1}^{-1}\ g_{2}^{-1}\ R_{12}(\lambda)\ g_{1}\ g_{2}. \label{solu3} \ee Use of the (assumed to exist) operator $\log g$ allows to explicitly solve (\ref{solu3}) as \be R_{12}(\lambda)
=\exp[-\sigma (\log g_{1} + \log g_{2})]\
(R_{12}^{0})\ \exp[\sigma (\log g_{1} + \log g_{2})] \label{solu4} \ee where again $\sigma$ denotes the sum over all dynamical variables $\sigma = \sum_{i=1}^n \lambda_{i}$ and $R^{0}$ does not depend on any variable $\lambda_{i}$ (except as usual, in dynamical Yang--Baxter equation, as an integer--period function). Note that in the example of \cite{ACF} where $g=\exp[{d\over du}]$ ($u$ spectral parameter), $R(\lambda)$ is again an exact adjoint action, so is (\ref{solu4}), hence $R^{0}$ is a $c$-number matrix.

Consider now (\ref{yben}a). From (\ref{solu4}) and (\ref{solu1}) one gets \be R_{12}^{0}\ R_{13}^{0\ gg}\ R_{23}^{0}= R_{23}^{0\ gg}\ R_{13}^{0}\ R_{12}^{0\ gg} \label{modyb} \ee hence $R^{0}$ is any non-dynamical solution of the shifted Yang--Baxter equation. It is in general not possible to go beyond this statement. However, particular solutions can easily be characterized. Any solution of the ordinary Yang--Baxter equation, commuting with $g \otimes g$ solves (\ref{modyb}). In the case described in \cite{ACF}, for instance $g = \exp[{d\over du}]$, any non-dynamical matrix
with a difference-dependance $R_{12}(u_{1} -u_{2})$ solves (\ref{modyb}).

To summarize, we now have the parametrized $A, B, C$ as \be && A_{12} = b_{1}^{-1}\ (g_{2}\ b_{2}\ g_{2}^{-1})^{-1}\ \Big \{Ad.\ \Big (\exp[-\sigma  (\log g_{1} +\log g_{2})]\Big )\ R_{12}^{0} \Big \}\ g_{1}\ b_{1}\ g_{1}^{-1} \label{newpa1} \\ &&
B_{12} =C_{21} = b_{2}^{-1}\ g_{2}\ b_{2}(\lambda +h_{1})\
g_{2}^{-1} \label{newpa2} \ee where $R^{0}$ solves
(\ref{modyb}). The existence of the $g$ shift in the Yang--Baxter equation (\ref{yben}b) coupled to the dynamical ``shift'' symbolized by ($h_{3}$) induces in the example in \cite{ACF} a 
coupling between the dependance in the dynamical parameters and the spectral parameter. Indeed (\ref{newpa1}) will read in this case: \be A_{12}(u_{1},\ u_{2},\ \lambda) =b_{1}(\lambda)\ b_{2}(u_{2} +\gamma, \lambda)\ R_{12}^{0}(u_{1} -\sigma,\ u_{2} - \sigma)\ b_{1}(u_{1} +\gamma)\ b_{2}(u_{2}). \label{acf}\ee

\subsection{The $D$-matrix and $K$ solutions}

A situation similar to section 3 arises here. Assuming first of all that $A,\ B,\ C$ are parametrized as in (\ref{newpa1})--(\ref{acf}) one is lead to discuss whether $D$ can be
\\
\\
1. detwisted at all or not
\\
\\
2. detwisted to a $g$--quasi non dynamical (QND) $R$-matrix
\\
\\
3. detwisted to a $g' \neq g$--QND $R$-matrix, where $g'$ may be simply $\1$, or any automorphism of ${\cal V}$. Situations 2 and 3 may once again overlap, but this problem will not be treated here.

One again establishes immediately that
\\
\\
{\bf Proposition 2'}
\\
If $ABC$ are parametrized as in (\ref{newpa1})--(\ref{acf}) by matrices $b$ and $R^0$ the following two
statements are equivalent:
\begin{itemize}

\item{The SDYB equation (\ref{SDRE}) has an invertible scalar solution $k(\lambda)$}
\item{$D$ can be detwisted according to (\ref{param02}) to a $g$-quasi non dynamical $R$-matrix $R$ with twist $q=gk$.}

\end{itemize}

\subsubsection{$D$ detwistable}

Let us first consider together cases 2 and 3 where $D$ can be rewritten as in (\ref{param02})
\be D_{12}= q_{1}^{-1}(\lambda+h_{2})\
q_{2}^{-1}(\lambda)\ \tilde R_{12}\ q_1(\lambda)\ q_{2}(\lambda+h_{1}) \label{dpar} \ee where
$\tilde R$ is a $g'$ quasi non-dynamical $R$ matrix i.e. obeys: \be
\tilde R_{12}(\lambda+h_{3}) = g_{1}^{'-1}\ g_{2}^{'-1}\ \tilde R_{12}\ g'_{1}\ g'_{2}. \label{req} \ee From (\ref{req}) it now follows that $\tilde R$ must obey the $g'$-modified non-dynamical Yang--Baxter equation
\be g_{1}^{'-1}\ g_{2}^{'-1}\ R_{12}\ g'_{1}\ g'_{2}\ R_{13}\ g_{2}^{'-1}\ g_{3}^{'-1}\ R_{23}\ g'_{2}\ g'_{3}= R_{23}\ g_{1}^{'-1}\ g_{3}^{'-1}\
R_{13}\ g'_{1}\ g'_{3}\ R_{12}. \label{gmod} \ee 
Eliminating now the non-trivial dynamical dependance of $R$ implied by (\ref{req}) we set \be \tilde R_{12} = \exp-\sigma (\log g'_{1} +\log
g'_{2})\ \bar R_{12}\ \exp\sigma (\log g'_{1} +\log g'_{2}) \label{eldy}\ee where now $\bar R_{12}(\lambda+h_{3}) = \bar R_{12}$, hence is independent (up to integer-periodic functions) of all dynamical variables. $\bar R$ also obeys the shifted non-dynamical Yang--Baxter equation (\ref{gmod}) equivalent ot (\ref{paramt2}).


Denoting now the quasi-non dynamical $R$-matrices respectively by
$R^{0\ d} \equiv $Ad exp $- \sigma \log g_1 + \log g_2 R^0$ (for $A$) and $\tilde R$ (for $D$) and plugging the corresponding parametrizations of $A, B, C, D$ into (\ref{SDRE}) one gets (denoting here by $K$ the solution of (\ref{SDRE}))
\be R_{12}^{0\ d}\ (g b g^{-1}
Kq^{-1})_{1}\ q_{1}\ (g b g^{-1} Kq^{-1})_{2}(h_{1})\
q_{1}^{-1} = (g b g^{-1} Kq^{-1})_{2}\ q_{2}\ (g bg^{-1}
K  q^{-1})_{1}(h_{2})\ q_{2}^{-1}\ \tilde R_{12}. \label{nsdr} \ee
As in (\ref{geq}) it is not easy to formulate general solutions $q$ to (\ref{nsdr}). However if the conjugations
$q_{i}X_{j}(h_{i}) q_{i}^{-1}$ can be trivialized, i.e. $(g b g^{-1} K q)_{i}(h_{j})$ is trivial on $V_{j}$, one can give explicit formulations of solutions $K$ in terms of non-dynamical objects ${\cal Q}$ by eliminating all dynamical dependence between $R^{0d}$ and $\tilde R$, rexpressing the equation in terms
of $R^0$ and $\bar R$. Consider first the case $g'=g$.

\subsubsection{$D$ detwistable to $g$-QND matrix}

As in the simpler case $g = \1$ two sets of solutions can be described:
\\
\\
{\bf Case 1. non-dynamical situation}
\\
\\
{\bf Proposition 4a}
\\
If ${\cal Q}^0$ is a non-dynamical solution to the non-dynamical shifted
Yang--Baxter equation: \be R_{12}^{0}\ {\cal Q}_{1}^0\ g_{2}^{-1}\
{\cal Q}^{0}_{2}\ g_{2} = {\cal Q}_{2}^0\ g_{1}^{-1}\ {\cal Q}_{1}^0\ g_{1}\ \bar
R_{12} \label{shift} \ee then \be K =g\ b^{-1}\ g^{-1}\Big
(\exp[-\sigma \log g]\  {\cal Q}^{0}\ \exp[\sigma \log g] \Big )\ q
\label{shso} \ee is also a solution of the SDR equation (\ref{SDRE}). If
$R^0 = \bar R$ there exists at least one invertible $ {\cal Q}^0 =\1$ and
$k = gb^{-1}g^{-1}q$ provides an invertible scalar solution to
(\ref{SDRE}) consistent with Proposition 2'.

More generally one has:
\\
\\
{\bf  Case 2: quasi non-dynamical solution}
\\
\\
Given an ad-factorizable automorphism $a$ on ${\cal V}$ such that $~[R^{0},\
a\otimes a]=[\bar R^{0},\ a \otimes a] =0$ one also establishes
\\
\\
{\bf Proposition 4b}
\\
If ${\cal Q}^{0}$ is a non-dynamical solution of the doubly shifted
RTT-type  equation \be R_{12}^{0}\ {\cal Q}_{1}^{0}\
(a_{2}^{\sigma} g_{2}a_{2}^{-\sigma})\ a_{2}^{-1}\ {\cal Q}_{2}^{0}\ a_{2}\ (a_{2}^{\sigma }g_{2}^{-1}a_{2}^{-\sigma }) ={\cal Q}_2^0\ (a_{1}^{\sigma }g_{1}a_{1}^{-\sigma })\ a_{1}^{-1}\ {\cal Q}_{1}^{0}\ a_{1}\ (a_{1}^{\sigma }g_{1}^{-1}a_{1}^{-\sigma}) \bar R_{12} \label{dsyb}\ee
(where $a^{\sigma} = \exp[\sigma\log a]$, assuming that $a$ has an operatorial logarithm) then \be K= g\ b^{-1}\ g^{-1}\ \exp[-\sigma \log g]\
\exp[\sigma \log a]\ {\cal Q}^{0}\ \exp[-\sigma  \log a]\ \exp[\sigma \log g]\ q(\lambda). \label{shso2} \ee is a solution of the SDRE (1.1).

In the particular case $\bar R = R^{0}$, (\ref{dsyb}) is
immediately solved by ${\cal Q}^{0} ={\1}$ hence (\ref{shso2})
defines an invertible solution to SDRE. If reciprocally one can identify an invertible solution ${\cal Q}^{0\ V}$ to (\ref{dsyb}), then $K$ provides an invertible solution to (1.1). As a consequence, as in the case of unshifted Yang--Baxter equations, $D$ can be directly detwisted to \be \tilde R_{12} = \Big ( Ad.\ \exp-\sigma (\log g_{1} 
+\log g_{2})\ R_{12}^{0} \Big )\ee using $(Ad.\ exp-(\sigma\ \log g)\ q)\ K$  as a twist instead of q in (\ref{dpar}).

\subsubsection{$D$ detwistable to a $f$-QND $R$-matrix, $f$ ad-factorizable}

The situation becomes here rather intricate. One can however show that extensions of
the two previous cases exist. Consider case 1. Equation (\ref{shift}) becomes
\be R_{12}^0\ {\cal Q}^{0}_{1}\ g_{2}^{-1}\ {\cal Q}_{2}^0\ f_{2} = {\cal Q}_{2}^{0}\ g_{1}^{-1}\ {\cal Q}_{1}^0\ f_{1}\ \bar R_{12}. \label{shift2b} \ee
Solutions are then given by:
\be K(\lambda) = g\ b^{-1}\ g^{-1}\ exp[-\sigma \log g]\ {\cal Q}^{0}\ exp[\sigma \log f]\ q(\lambda). \label{sol2b}\ee
Case 2 can be also extended to this case. The relevant equations become: \be R_{12}^0\ {\cal Q}_{1}^{0} \Big \{Ad. exp[-\sigma \log a]\ g \Big \}^{-1}_2\ \Big ( a^{-1}{\cal Q}^{0}a \Big )_{2}\ \Big \{Ad. exp[-\sigma \log a]\ f \Big \}_2 = \nonumber\\ {\cal Q}_{2}^0\ \Big \{ Ad. exp[-\sigma \log a]\ g \Big \}^{-1}_{1}\ \Big ( a^{-1}{\cal Q}^{0}a \Big )_{1}\ \Big \{ Ad. exp[-\sigma \log a]\ f \Big \}_{1}\ \bar R_{12}. \label{ybl} \ee and  solutions are given by:
\be K(\lambda) = g\ b^{-1}(\lambda)\ g^{-1}\ exp[-\sigma \log g]\ exp[-\sigma \log a]\ {\cal Q}^0\ exp[\sigma \log a]\     exp[\sigma \log f] q(\lambda). \label{sol3b} \ee
\\
\\We must make two important remarks here:

First of all it is important to notice that in both equations (\ref{ybl})
and (\ref{dsyb}) an explicit conjugation of the ${\cal V}$--automorphism $g$ by a dynamical ${\cal V}$--automorphism
$exp[\sigma \log a]$ occurs. If $[a,\ g]=0$ no conjugation occurs and (\ref{ybl}), (\ref{dsyb}) are 
genuine non-dynamical Yang--Baxter RTT type equations for which it is consistent to search for
non--dynamical solutions ${\cal Q}^0$. If not it may be impossible to find non-dynamical solutions ${\cal Q}^0$ 
and these cases may then be empty.

Second remark: Once again if $exp[\sigma \log g]$ or $exp[\sigma \log f]$ are not factorizable, even though
$f$ and $g$ are ad-factorizable, the RTT-type equations are not written as finite-size matrix
algebraic equations on tensor products of the auxiliary space $V$. Solutions ${\cal Q}^0$
may then not be finite-size matrices and may not admit an expansion as formal
power series of the variable $u$; and the object $K$ may not be viewed as generating functional
of some quantum reflection-like algebra.

\subsubsection{$D$ not detwistable}

If no parametrization of $D$ can be defined on the lines of  (\ref{param02}), one can again still prove
the comodule property:
\\
\\
{\bf Proposition 5}
\\
If $K(\lambda)$ is a solution of (1.1) and ${\cal Q}$ is a non-dynamical
solution of \be R_{12}^{0}\  {\cal Q}_{1}\ (g^{-1}\  {\cal Q}_{2}\ g) = {\cal Q}_{2}\
(g^{-1}\  {\cal Q}_{1}\ g)\ R_{12}^{0} \label{qs}\ee where $R^{0}$
and $g$ are defined in (\ref{newpa1})--(\ref{newpa2}), then \be (g_{1}\ b_{1}^{-1}\
g_{1}^{-1})\ \Big (\exp[-\sigma \log g]\  {\cal Q}_{1}\
\exp[\sigma \log g]\Big )\ g_{1}\ b_{1}\ g_{1}^{-1}\ K \label{ns}\ee  is also a solution of (1.1). The dressing of an a priori given (operatorial or scalar) solution$K(\lambda)$ by suitable ``dynamical'' solutions of the Yang--Baxter equation (\ref{qs}) seems to be the only available construction of new solutions in this case.

We shall now give explicit simplified formulations for
the monodromy matrices obtained from the simplest comodule
structures defined in \cite{NAR}, in the simplest parametrization
context defined by Proposition 4a.

\subsection{Monodromy matrices when $D$ detwistable to $g$-QND $R$}

When $D$ can be detwisted to a quasi-non dynamical $R$ of the same type as $A$ the monodromy matrix built by using the comodule structure of the SDYB reflection equation, with appropriate $\tilde A, \tilde B, \tilde C, D$ matrices and a scalar solution $k(\lambda)$, will again simplify. Let us first consider the simplest case where $D$ is detwisted to the same matrix $\tilde R$ as $A$, equivalent to the existence of invertible scalar solutions to (\ref{SDRE}). One defines the consistent parametrization: 
\be 
&& A_{12} = b_{1}^{-1}\ (g_{2}\ b_{2}\ g_{2}^{-1})^{-1}\ \tilde A_{12}\ b_{2}\ (g_{1}\ b_{1}\ g_{1}^{-1})
 \label{3par1} \\ 
&& \tilde A_{12} =Ad.\ \exp[-\sigma  (\log g_{1} +\log g_{2})]\ R_{12}^{0} 
\label{3par2} \\ 
&& B_{12}= C_{12}^{\pi}= b_{2}^{-1}\ g_{2}\ b_{2}(h_{1})\ g_{2}^{-1}
\label{3par3} \\ 
&& D_{12} = k_{1}^{-1}(h_{2})\ C_{12}^{-1}\ k_{2}^{-1}\ A_{12}\  k_{1}\ B_{12}\ k_{2}(h_{1}). \label{3par4} 
\ee 
Eq. (\ref{3par4}) just reflects the fact that since $D$ is detwisted to $\tilde R =\tilde A$, there exists an invertible scalar solution $k$ to (\ref{SDRE}) which can be used directly to rewrite $D$.

The monodromy matrix for a $N$-site chain is now defined once one stipulates a direct $ {\cal Q}_0$ and a dual $\chi_{0}$ (scalar) reflection matrix. We choose for $ {\cal Q}_0$ the simplest parametrization
described by Proposition 4a when $R$ and $\bar R$ are identical. To define the dual
solution $\chi_{0}$ we use the known identification between transposed solutions of dual SDRA,
and inverse of direct solutions of SDRA, meaningful here since we know from Prop. 2' that such invertible
solutions exist. 
We set accordingly:

\be &&  {\cal Q}_{0} =g_{0}\ b_{0}^{-1}\ g_{0}^{-1}\ \tilde  {\cal Q}_{0}\ g_{0}\ b_{0}^{-1}\ g_{0}^{-1}\ k \label{set1} \\ && \tilde  {\cal Q} = Ad.\ \exp[-\sigma  \log g]\  {\cal Q}_{R}^{0}, ~~~~\mbox{where} ~~~~  {\cal Q}_{R}^{0} ~~~\mbox{obeys (\ref{shift})} \label{set2} \\ && \chi_{0}^t = k_{0}^{-1}\ g_{0}\ b_{0}^{-1}\ g_{0}^{-1}\ \tilde  {\cal Q}'_{0}\ g_{0}\ b_{0}^{-1}\ g_{0}^{-1} \label{set3} \\ && \tilde {\cal Q}' = Ad.\ \exp[-\sigma \log g]\  {\cal Q}_{L}^{0\ -1} ~~~~\mbox{where} ~~~~  {\cal Q}_{L}^{0} ~~~\mbox{obeys (\ref{shift})}. \label{set4} \ee

The monodromy matrix now reads \cite{NA}: \be {\cal T}_{0}\
e^{\partial_{0}} &=&  \chi_{0}^t\ g_{0}\ A_{0\ 2N}\ g_{2N}^{-1}\ g_{0}\ C_{0\ 2N-1}\ g_{0}\ A_{0\ 2N-2}(h_{2N-1})\ g_{2N-2}^{-1} \ldots \\ && \ldots  {\cal Q}_{0}(h_{1}+h_{2} +...h_{2N-1})\ D_{01}(h_{1}+h_{2}
+...h_{2N-1}) \ldots g_{2N}\ B_{0\ 2N}\ e^{\partial_{0}} \label{mono3}\ee {\bf Remark:} Contrary to the scalar $k$ matrix, the matrix ${\cal T}_{0}\ e^{\partial_{0}}$ exhibits a non-adjoint action of $g_{0}$ (but an adjoint action of all non-zero indexed operators $g_i$). This may lead to a fundamental problem:

In the non-affine case, when $dim {\cal V} \equiv V < \infty$ the transfer matrix is defined as a trace over $V$ hence no difficulty arises. If however ${\cal V}$ is an evaluation module $V \otimes {\mathbb C}[[u]]$, one is actually interested in partial traces over $V$ to define spectral-parameter dependent transfer matrices $Tr_{V}({\cal T}_{0}\ e^{\partial_{0}})$. In this case if $g$ acts non-trivially on ${\mathbb C}[[u]]$ , more specifically
if $g$ is not factorizable, (as in \cite{ACF}  where
$g=\exp[{d\over du}]$) the proof of commutation of such partial traces using the $A{\cal T}B{\cal T}$ relations is not valid, as can be seen on our example since the ${\cal T}$ matrices will then contain explicit operators 
$\exp[{d\over du}]$ acting on matrix elements of $A, B, C, D$! As a matter of fact even the partial traces over such monodromy matrices do not
exist since the matrices themselves do not assume the factorized form of dim ($(V)^{\otimes N} \otimes V$)-size matrices 
depending on $N+1$ spectral parameters. 

A solution to this issue is the following: One has to assume that $D, B, C$
exhibit the same zero-weight properties under the adjoint action of $g$ as they already did, 
as a fundamental assumption of our semi-dynamical structure,
under the adjoint action of ${\mathfrak h}$. In addition one will assume that $g$ and ${\mathfrak h}$ commute: \be \Big [D,\ g \otimes g  \Big ]= \Big [B_{12},\
g \otimes {\1} \Big ] = \Big [C_{12},\  {\1} \otimes g \Big ] =0 ~~~~~\Big [h,\ g \Big] =0. \label{zwc}\ee 

This situation is indeed realized
in \cite{ACF}  since (\ref{zwc}) here immediately follows from the particular
dependence of $D,\ B$ and $C$ on the spectral
parameter: $g=\exp[{d\over du}]$ and $D_{12} =D_{12}(u_{1} -u_2)$,
$B_{12} =B_{12}(u_{2})$ and $C_{12} = C_{12}(u_1)$. Adjoint action
of $g$ is simply shift of the corresponding spectral parameter.
 
Once (\ref{zwc}) is imposed it is easy to prove:
\\
\\
{\bf Proposition 6}
\\
If $K$ is a solution to (\ref{SDRE}), $K\ g^{n}$ is a solution to (\ref{SDRE}) for any integer $n \in {\mathbb Z}$.
\\
\\
The monodromy matrix (\ref{mono3}) can then be modified to take the form 
of an exact adjoint action (hence factorizable) : \be {\cal T}_{0}\ e^{\partial_0} \to {\cal T}_{0}\
e^{\partial_0}\ g_{0}^{-2N}. \label{norm}  \ee Since we have
restricted $g$ to be ad-factorizable, the partial trace of the monodromy matrix is now
once again correctly defined; its expansion in formal power series of $u$ is also defined and  
generates commuting Hamiltonians. 

Let us make here a technical remark: locality conditions on
these Hamiltonians may then be imposed (see \cite{NA}) and lead to specific
choices of the values of the quantum--space spectral parameters:
As a particular example let us point out that in the case treated
in \cite{ACF}, 
the shifts in (\ref{mono3}) are distributed 
according to: \be  && \ldots A_{0\ 2n}(\lambda_{0} +
(1+2N-2n), \lambda_{2n})\ C_{0\ 2n-1}(\lambda_{0} +
(2+2N-2n), \lambda_{2n-1}) \nn\\ && \ldots D_{0\ 2n-1}(\lambda_{0}
+(2N), \lambda_{2n-1})\ B_{0\ 2n}(\lambda_{2n} )
\label{distr} \ee and the locality conditions on the Hamiltonians
have a consistent implementation as $\lambda_{2n} =\lambda_{0}
+(2N -2n+1),\ \lambda_{2n-1} = \lambda_{0} +2N$.

If these assumptions are realized, plugging now (\ref{3par1}),
(\ref{3par3}), (\ref{set1}), (\ref{set3})  into (\ref{mono3})
yields \be {\cal T}_{0} e^{\partial_{0}}  = {\cal O}_{N}^{-1}\ \tilde {\cal T}_{0}\
e^{\partial_{0}}\ {\cal O}_{N} \label{mono4} \ee 
where \be \tilde
{\cal T}_{0} &=& k_{0}^{-1}\ g_{0}\ b_{0}^{-1}\ g_{0}^{-1}\ \tilde
 {\cal Q}'_{0}\ g_{0}\ \tilde A_{0\ 2N}\ \ldots  g_{0}\  A_{02}\
\tilde  {\cal Q}_{0}\ g_{0} \tilde A_{01} \ldots g_{0} \tilde A_{0\ 2N-1}
\nn\\ && (g_{0}\ b_{0}\ g_{0}^{-1})\ k_{0}\ g_{0}^{-2N}
\label{tilde1}\ee 

\be {\cal O}_{N}= \Pi_{m=0}^{N-1}(g_{2N-2m}\ b_{2N-2m}(h_<^{odd})\ g_{2N-2m}^{-1})\
(g_{2N-2m-1}\ b_{2N-2m-1}(h_<^{odd})\ g_{2N-2m-1}^{-1}\
k_{2N-2m-1}(h_<^{odd})) \nn\\ \label{calo}\ee 

Reformulating the quasi-non
dynamical $\tilde A$, $\tilde Q_{0}$ and $\tilde Q_{0}'$ in (\ref{tilde1}) following
(\ref{3par1})--(\ref{set4}) one finally gets
 
 \be \tilde {\cal T}_{0} &=&
k_{0}^{-1}\ g_{0}\ b_{0}^{-1}\ g_{0}^{-1}\ \exp[-\sigma  \log
g_{0}]\  {\cal Q}_{L}^{0\ -1} Ad \exp[-\sigma (\log g_{1}+...g_{2N} )] \nn\\
&& \{ g_{0}\ R^{0}_{0\ 2N} \ldots g_{0}\ R^{0}_{02}\  {\cal Q}_{R}^{0}\ g_{0}\ R^{0}_{01}
\ldots g_{0}\ R^{0}_{0\ 2N-1} \} \exp[\sigma  \log g_{0}]\ g_{0}\ b_{0}\
g_{0}^{-1}\ k_{0}\ g_{0}^{-2N}. \label{tilde2} \ee {\bf Comment:}
$\tilde {\cal T}_{0}$ is therefore decomposed as a non-dynamical
chain monodromy matrix with direct/dual ``reflection'' matrix
dressed dynamically by the shift-dynamical coupling $Ad.\ \exp
[-\sigma \log g_{0}]$; more fundamentally dressed by the adjoint
action of the Drinfeld twist $g_{0}\ b_{0}\ g_{0}^{-1}\ k_{0}$,
which turns $D$ into $\tilde R$), yielding a generating
functional for the commuting Hamiltonians by the dynamical
trace formula  $~Tr_{0} \Big (\tilde
{\cal T}_{0}\ \exp[\partial_{0}] \Big)$.

\subsection{Monodromy matrices when $D$ detwistable to $\bar R$ not equivalent to $R$}

Here one must substitute to (\ref{3par4}) the general twisting relation (\ref{dpar}). Using now as direct reflection matrix a solution of the form (\ref{shso}) and a (non parametrized) dual reflection matrix $\chi_{0}^t$ one gets for the monodromy matrix a formula analogous to (\ref{tilde2}) with the following modifications:

\begin{itemize}

\item{1. The blocks $(g b g^{-1} k)$ in ${\cal O}$ 
and the term $g_{0}\ b_{0}\ g_{0}^{-1}\ k_{0}$
on the r.h.s. of (\ref{tilde2}) must be 
substituted by the twist matrix $q$ from $D$ to $R$.}

\item{2. Odd-labelled $R_{0\ 2k+1}$ are substituted by $\bar R_{0\ 2k+1}$ defined in (\ref{eldy}).}

\item{3. Since no invertible solution to  (\ref{shift})
exists here, we cannot
identify a dual solution with any ``inverse'' of a direct solution. 
Parametrization (\ref{set3}, \ref{set4}) is however still valid provided
that $k_{0}^{-1}\ g_{0}\ b_{0}^{-1}\ g_{0}^{-1}$ be replaced by $q_0^{-1}$
and ${\cal Q}_{L}^{0\ -1}$ by an explicitely computed solution
of the transposed dual equation to (\ref{shift}). This transposed dual equation
is trivially obtained by taking the formal
inverse of (\ref{shift}). The term $k_{0}^{-1}\ g_{0}\ b_{0}^{-1}\ g_{0}^{-1}$
on the l.h.s. of (\ref{tilde2}) is consequently to be substituted by $q_0^{-1}$.}

\end{itemize}
It is not clear whether such transfer matrices are useful to build
physically interesting spin chain type models. Their explicit
formulation however may be interesting in itself to understand
the algebraic structures underlying (1.1) in the non-trivial case
where $A, B, C$ and respectively $D$ yield distinct $R$ matrices.

\subsection{Remarks on the structure of monodromy matrices}

As commented upon in the previous sections, the monodromy matrices
take a very characteristic form once the parametrization of $A, B,
C, D$, $k$ and $k^{dual}$ is taken into account. One identifies
first non-dynamical chain transfer matrices with direct and dual
scalar Lax matrix $q_{R}$ and $q_{L}$, which would yield by the
standard construction closed spin chain Hamiltonians. They are
then dressed non-trivially by the adjoint action of the Drinfeld
twist $q$, which characterizes the $D$ matrix, and the subsequent
generating functional for commuting Hamiltonians is: \be
t(\lambda)= Tr_{0}\{ q_{0}^{-1}(\lambda)\ T_{0}\ q_{0}(\lambda)\
e^{\partial_{0}} \} \label{transfer} \ee The key remark here is
that mutual commutation of such objects with different spectral
parameters $u_{0},\ u_{0}'$; or of ``quantum trace-like'' objects
obtained from fusion procedures on the auxiliary space ((0)
-index) as was derived in \cite{NADR}; is guaranteed by the necessary
conditions on the twist $q_{0}$, i.e. that the $D$-matrix obtained
as dynamical twist of the non-dynamical or quasi-non dynamical
$R$-matrix in ${\cal T}_{0}$ as: \be D_{12} = q_{2}^{-1}(h_{1})\
q_{1}\ R_{12}\ q_{2}\ q_{1}(h_{2}) \label{dec} \ee have zero
weight. Otherwise the ($\exp[\partial_{0}]$) term prevents
commutation of the generating functions ($[t(u),\ t(u')] =0$).
Remarkably though, zero weight condition on $D$ is also a
sufficient condition (Proposition 1) to guarantee that $D$
obeys the dynamical Yang--Baxter equation.

This leads us to conclude that the semi-dynamical ``reflection''
equation is not really a ``reflection'' equation, in the usual
sense of the term, since in any case $B$ and $C$ have
non-canonical, loosely speaking ``semi-diagonal'',
zero-weight conditions. It seems that one underlying
fundamental structure is the dynamical Yang--Baxter algebra
(dynamical quantum group) associated to the matrix $D$; the
decomposition (\ref{dec}) is then used to build dynamical
monodromy matrices (\ref{transfer}) although bypassing the
zero-quantum weight requirement \cite{ABB}, which occurs when using
directly Lax matrices of the dynamical quantum group to build monodromy matrices. This
requirement is eliminated by the trick of building a
reflection-type quadratic exchange algebra with no dynamical
shifts in the coefficient matrices. The SDRE is therefore an
intermediate construction between the non-dynamical quantum group
($R$-matrix) and the dynamical quantum group ($D$-matrix). Its
main practical interest is that it naturally yields a dynamical
un-constrained (no-zero weight) monodromy matrix (\ref{transfer}).
Let us once again remark that in any case (\ref{transfer}) does
not admit an obvious interpretation as a trace in the quantum
group. The groupoid formulation, advocated in \cite{Ping}, and quite
naturally adapted to dynamical $R$ matrices may provide a natural
framework for (\ref{transfer}).

\section{Conclusions and perspectives}

\subsection{The spin-chain Hamiltonians}
Traces over the auxiliary space (labelled by $0$) of the
monodromy matrices such as described in Sections
3.3 and 4.3 provide a systematic way of constructing
quantum integrable Hamiltonians \cite{NA}. 

It is first of all essential to remark
that in this construction
the quantum-space operators ${\cal O}_{N}$
are not relevant to keep
since they simply conjugate the quantum monodromy matrix,
and the Hamiltonians deduced from it.
One must therefore realize
the computation of quantum
integrable Hamiltonians from the 
non-conjugated monodromy matrix, thereby eliminating
all cumbersome quantum-space shifts. 

These now factored--out monodromy matrices exhibit
a very interesting combination of features. The untwisted part 
$R ...{\cal Q} ...R$, as already mentioned,
has the canonical form of a generating
functional (once taking the trace over the
auxiliary space) for closed spin--chain Hamiltonians. The
zero--site twisted monodromy matrix, built
from a single scalar reflection matrix $k$ and a trivial
dual solution $\1$, yields
precisely scalar RS Hamiltonians
when choosing $ABCD$ structure matrices from \cite{ACF}.
The question of how such features interplay in the new generated
Hamiltonians to yield possible ``spin Ruijsenaar Schneider models''
is therefore quite challenging. More precisely the procedure 
should now run as follows.

In the non--affine case (no spectral parameter) the trace over
the full auxiliary space $V$ 
is expected to yield Hamiltonians of $N$-body systems
with interactions a la Ruijsenaar--Schneider. 
A family of commuting higher--degree Hamiltonians
can then be obtained by a now well-established quantum
trace procedure, see e.g.\cite{ABB,NADR}. 

Consider now the more interesting case
of affine SDRA where ${\cal V} = V\otimes {\mathbb C}[z]$. 
It was shown in \cite{NA} how
quantum integrable ``spin RS'' Hamiltonians,
could be obtained by the canonical procedure of 
taking the logarithmic derivative
of the partial trace over $V$ of the dynamical
monodromy matrix, 
w.r.t. the spectral parameter $z_0$ associated to
the auxiliary space, at $z_0 = 0$. 

The generalized ``spin-spin'' interactions then
take a local form (nearest neighbor or next-to-nearest neighbor
interaction) for a suitable consistent choice of the values
of the ``quantum'' spectral parameters $z_{1...2N}$,
provided that the structure matrices
$A$ and $D$ obey the so-called ``regularity'' conditions
$A(z_1 = 0, z_2 = 0) = D(z_1 = 0, z_2 = 0) = P_{12}$,
where $P_{12}$ is the permutation operator on $V\otimes V$.
This suitable choice was commented on in Section 4, see \ref{distr}.

The technical problem which arose in the previous approach \cite{NA} 
when directly computing the form of these Hamiltonians lied
in the complexity of the formulae once written in terms of non--parametrized
matrices $A, B, C, D$. 
Reformulating the structure matrices 
as we have done, the new
factorized monodromy matrices are essentially formulated in terms
of one single non--dynamical $R$--matrix $R$
and one consistently associated twist matrix $q$ yielding
a dynamical $D$ matrix through \ref{paramt}. The affine case
is however the situation where \cite{BRT} in principle does not extend, and $q$-matrices
must be explicitly constructed ``by hand''
from given $D$-matrices. It is a priori
known that they exist for the specific
RS $ABCD$ matrices
since in this case all $d_{ij}$ are non zero, all $b_i$
are invertible, and $k(\lambda) = \1$ is known to be a
solution. Proposition 2' then applies; the
final step to get explicit Hamiltonians
is now to compute matrices $q$ and $R$
from known non-constant $D$ matrices
and derive explicit ``spin-chain
RS''-type Hamiltonians from the factorized
forms (\ref{tilde2}).  

The specific form of the interaction will also
depend on the choice of the scalar solution $k$. 
Classification of solutions $k$ for the non-affine rational
case of \cite{ACF} is now fully known \cite{AR}, and classification
in the affine rational case is currently in progress.

We finally want to indicate that relaxing some technical
hypothesis, such as non-factorizability, may
yield interesting generalizations of the RS
type Hamiltonians \footnote{This intriguing possibility was suggested
to us by one referee}.

\subsection{Connection between SDRA and quantum groups}

In the most regular situations, when $B =C^{\pi}$ is invertible,
$d_{ij} \neq 0$ for all $i,\ j$ and SDRE (\ref{SDRE}) has at least one
invertible scalar solution $k(\lambda)$, the simplest parametrization 
of solutions $K$ proposed as subset $3.2.2 a1$
allows to prove that any representation of the ordinary quantum
group ($RTT=TTR$) generates a representation of the SDRA. The
inference is only one--sided since one cannot preclude the
possibility of dynamical solutions to the reduced equation (\ref{geq}). 
In addition the monodromy matrices generated by the simplest
representations of the comodule structure by $ABCD$ matrices are 
expressed in terms of the standard monodromy matrix generated by
the ordinary quantum group $R$-matrix, which yields closed
spin-chain Hamiltonians. In this respect
one can say that the twisting procedure is compatible with the comodule
structure, which is of course to be expected if they 
represent universal algebraic structures (Drinfel'd twist
and coproduct) which would underlie the SDRA. 
 
 The form of the dynamical trace $Tr\{{\cal T}
e^{\partial_0} \}$ however remains a specific feature of the SDRE,
and --as already commented-- does not naturally yield an element
of the SDRA itself, but may rather be understood in terms
of a more complex algebraic structure, possibly a quantum groupoid \cite{Ping,Et}.

These conclusions can be extended to the ``$g$-modified''
extension of the SD Yang--Baxter equations. The ``$g$-modified
quantum group'' structure $R_{12}T_1T_2^g= T_2 T_1^g R_{12}$
however is a less standard one and certainly deserves more
exploration, in particular since it is the relevant one
when considering the elliptic Ruijsenaars-Schneider example developed
in \cite{ACF}. 

In the
non-regular situation, when no invertible $k(\lambda)$ is
available, the representations of intertwining relations
$R_{12}T_1T_2^g= T_2 T_1^g R_{12}$ are now relevant to build
representations of the corresponding SDRA, and monodromy matrices.
In fact, as mentioned before, the SD``R''E (\ref{SDRE}) is not so much defining
a reflection algebra as providing an intertwining formulation
between a conjugated $R$--matrix $A$ and a dynamical twisted $D$ matrix
with same or different underlying non--dynamical or quasi--non--dynamical
$R$--matrices, themselves associated to quantum group-like
algebraic structures.
A better understanding of this structure
may require a (partial) lifting of the sufficient conditions, e.g.
the (quasi) non-dynamicity condition on $gKg^{-1}$. In addition,
lifting the conditions of $B$-invertibility or
$d_{ij} \neq 0$ may provide interesting non-trivial new examples.

\appendix

\section{Appendix: Semi-dynamical quantum reflection algebra}

Quantum reflection algebras were first formulated in \cite{Che,Skl} as
consistency conditions between factorizable 2-body $S$-matrices of
quantum integrable systems, and 1-body reflections $K$-matrix,
guaranteeing the quantum integrability of the system with
boundaries. They take the general form \be A_{12}\ K_1\ B_{12}\
K_2 = K_2\ C_{12}\ K_1\ D_{12}.  \label{11} \ee Equations
(\ref{11}) is now interpreted as quadratic constraint equation for
generators of the quantum algebra ${\cal G}$ encapsulated in the
matrix $K$. It is represented as an equation in $\mbox{End} ({\cal
V}) \otimes \mbox{End} ({\cal V})$ with elements in $U({\cal G})$
where ${\cal V}$ is a given vector space known as the auxiliary
space. $ABCD$ are matrices in $\mbox{End} ({\cal V}) \otimes
\mbox{End} ({\cal V})$. ${\cal V}$ may be --in the most usual
case-- a finite vector space $V$ or a loop vector space $V \otimes
{\mathbb C}[[z]]$; (the abstract formal variable $z$ being the so
called ``spectral parameter''). However one may retain the
possibility that ${\cal V}$ be a more general vector space
(functional space), even though it will not be considered
in the present work. $K$ now belongs to $\mbox{End} ({\cal V}) \otimes {\cal
G}$. A generalized quantum reflection algebra may be defined when
assuming that $ABCD$ and $K$ depend on a further set of complex
variables, collectively denoted $\lambda = \{\lambda_i, ~~i=1, \ldots n \}$,
interpreted as coordinates on the dual of a characteristic
(usually abelian) complex Lie algebra of finite dimension, and
parametrizing a deformation of (\ref{11}). This is in fact an
extension to (\ref{11}) of the so-called dynamical deformation of
YB equation defined in \cite{Fe,GN,ABRR,JKOS} 
where the YB equation originally introduced as:
\be &&R_{12}\ R_{13}\ R_{23}=
R_{23}\ R_{13}\ R_{12} \\ &&R_{12}(\lambda + h_3)\ R_{13}\
R_{23}(\lambda+h_1)= R_{23}\ R_{13}(\lambda +h_2)\ R_{12}.
\label{22} \ee Here since $\lambda$ are coordinates on the dual of
${\mathfrak h}$, it is understood that the auxiliary space ${\cal
V}$ is an irreducible diagonalizable module of ${\mathfrak h}$,
justifying the notation ``$\lambda + h_{i}$''. ``Irreducibility''
is an extra requirement, implying that zero-weight matrices under
adjoint action of ${\mathfrak h}$ necessarily admit an expansion
of the finite generators of ${\mathfrak h}$, which will be very
useful in all our discussions.

In fact two dynamical extensions of the RA (\ref{11}) have now
been identified. The semi-dynamical RA, which interests us here,
reads: \be A_{12}\ K_1\ B_{12}\ K_{2}(\lambda+h_1) = K_2\ C_{12}\
K_{1}(\lambda+h_2)\ D_{12}. \label{12} \ee The fully dynamical RA
or ``boundary dynamical RA'' \cite{FHS,FHLS} reads \be  A_{12}\
K_1(\lambda+h_2)\ B_{12}\ K_{2}(\lambda+h_1) = K_2(\lambda+h_1)\
C_{12}\ K_{1}(\lambda+h_2)\ D_{12}. \label{13} \ee


\section{Remark: The irreducibility criterion in the affine case}

We have chosen two specific cases for the auxiliary space ${\cal
V}$, either as a finite dimensional vector space ${\cal V} = V$,
or as loop space ${\cal V} = V\otimes {\mathbb C}[z]$. The finite
dimensional vector space $V$ is in addition assumed to be a
diagonalizable irreducible module for the dynamical Lie algebra
${\mathfrak h}$, allowing in this way to consistently expand any
zero weight matrices on any basis of ${\mathfrak h}$, e.g. the
basis of normalized diagonal $n$-matrices $\{e_{ii} \}$ (here
$(e_{ij})_{kl} = \delta_{ik} \delta_{jl}$) when ${\mathfrak h}
=Cartan(gl(n))$.

It would seem that therefore, when ${\cal V} = V\otimes {\mathbb
C}[z]$, the full auxiliary space ${\cal V}$ is no more an irreducible
module of ${\mathfrak h}$. However if ${\mathfrak h}$ is completed
by the derivation generator $d$ ---as it should be when
considering affine Lie algebras--- represented as $d =
{d\over d u}$, ${\cal V}$ is again irreducible. One would expect
in this case occurrence of an $(n+1)$-th coordinate
$\lambda_{d}$ with a dynamical shift in (\ref{SDRE}). However, in the
known case of dynamical elliptic quantum groups \cite{JKOS}
the dynamical shift on the coordinate associated to
$d$ is interpreted as the central charge $c$ in a centrally
extended dynamical quantum algebra, hence it is set to 0 in an
evaluation representation. Since the shifts in the definition of
the dynamical reflection algebra (\ref{SDRE}) occur precisely on the auxiliary
spaces, the absence of an explicit $(n+1)$-th shift in (\ref{SDRE}) does not 
contradict the existence of a (here non relevant) extra variable (such as the 
elliptic module $p$ in an elliptic DRA) and thus the interpretation of (\ref{SDRE}) as
dynamical quantum reflection algebra, with dynamical Lie algebra
${\mathfrak h} \cup \{d \}$, for which ${\cal V}$ is again an irreducible
module. Note that the choice of 
$\hat{\mathfrak h}= {\mathfrak h} \cup \{ d\}$ as underlying
abelian Lie algebra defining the dynamical deformation now implies
---for consistency of the construction--- to implement 
full $\hat{\mathfrak h}$ zero-weight conditions on $B, C, D$, i.e.
including the adjoint action of $d$. In this case $g=d$ becomes a suitable automorphism,
under the conditions in (\ref{zwc}) to build
monodromy matrices in the $g$-deformed YB frame. This is precisely the
situation realized in \cite{ACF}.
\\
\\
{\bf Acknowledgements}

We wish to thank warmly A. Doikou for discussions and
overall help in the formulation of this paper. We thank
the referees for their helpful suggestions. J.A.
thanks INFN Bologna and Francesco Ravanini for
their hospitality. G.R. acknowledges ANR contract
JC05-52749 for partial support.

\end{document}